\documentclass[12pt]{article}

\usepackage{amsmath,amsthm,amsfonts,amssymb,amscd}

\usepackage[english,french]{babel}

\usepackage{epsfig}
\usepackage{graphicx}
\usepackage{multirow}
\usepackage{array}
\usepackage[all]{xypic}

\newtheorem{theorem}{Theorem}[section]

\newtheorem{corollary}[theorem]{Corollary}

\def\ml{\mathcal{C}}
\def\ml1{\mathcal{C}^1}
\def\mlb1{\mathcal{C}_{b}^{1}}

\def\mc{\mathbb{C}}

\def\PP{\mathbb{P}}

\usepackage{color}
\usepackage{ifthen}

\newboolean{usecolor}
\setboolean{usecolor}{true}

\ifthenelse{\boolean{usecolor}}
{
  \definecolor{colore}{cmyk}{0,1,0.6,0}
  \definecolor{coloregen}{cmyk}{0.7,0,1,0}
  \definecolor{coloresimo}{cmyk}{1,0.6,0,0}
}
{
  \definecolor{colore}{cmyk}{0,0,0,1}
  \definecolor{coloregen}{cmyk}{0,0,0,1}
  \definecolor{coloresimo}{cmyk}{0,0,0,1}
}

\title{On the Configuration Spaces of Grassmannian Manifolds}

\author{Sandro {\sc Manfredini}\footnote{Department of Mathematics, University of Pisa, manfredi@dm.unipi.it}  \and Simona {\sc Settepanella}\footnote{LEM, Scuola Superiore Sant'Anna, Pisa, s.settepanella@sssup.it. }}
\begin{document}

\maketitle

\begin{abstract}Soit $\mathcal{F}_h^i(k, n) $ le $i$\`eme espace de configuration ordonn\'e  de tous les points distincts $H_1,\ldots,H_h$ dans la Grassmannienne $Gr(k, n) $ de sous-espaces de dimension $k$ de $\mc^n$, dont la somme est un sous-espace de dimension $i$. Nous prouvons que $\mathcal{F}_h^i(k, n)$ est (si non vide) une sous-vari\'et\'e complexe de $Gr(k, n)^h$ de dimension $i(n-i)+hk(i-k)$ et son groupe fondamental est trivial si $i=min(n,hk)$, $hk \neq n$ et $n> 2$ et \'egal au groupe de tresses de la sph\`ere $\mc P^1$ si $n= 2$. Finalement, nous calculons le groupe fondamental dans le cas particulier des arrangements d'hyperplans, c'est \`a dire $k=n-1$.
\end{abstract}

\selectlanguage{english}

\begin{abstract}
Let  $\mathcal{F}_h^i(k,n)$ be the $i$th ordered configuration space of all distinct points $H_1,\ldots,H_h$ in the Grassmannian $Gr(k,n)$ of $k$-dimensional subspaces of $\mc^n$, whose sum is a subspace of dimension $i$. We prove that $\mathcal{F}_h^i(k,n)$ is (when non empty) a complex sub\-ma\-ni\-fold of $Gr(k,n)^h$ of dimension $i(n-i)+hk(i-k)$ and its fundamental group is trivial if $i=min(n,hk)$, $hk \neq n$ and $n>2$ and equal to the braid group of the sphere $\mc P^1$ if $n=2$. Eventually we compute the fundamental group in the special case of hyperplane arrangements, i.e. $k=n-1$.\end{abstract}

\begin{center}
{\small\noindent{\bf Keywords}:\\
complex space, configuration spaces, \\braid groups.}
\end{center}

\begin{center}
\end{center}

\section{Introduction}

Let $M$ be a manifold. The \emph{ordered configuration space} $$\mathcal{F}_h(M)=\{(x_1,\ldots,x_h)\in
M^h|x_i\neq x_j,\,\,i\neq j\}$$ of $h$ distinct points in $M$ has been widely studied after it has been introduced by Fadell and Neuwirth \cite{FN} and Fadell \cite{F} in the sixties. It is well known that for a simply connected manifold $M$ of dimension greater or equal than $ 3$, the \emph{pure braid group} $\pi_1(\mathcal{F}_h(M))$ on $h$ strings of $M$ is trivial. This is not the case when the dimension of $M$ is lower than $3$ as, for example, the pure braid group of the sphere $S^2\approx\mc P^1$ with presentation: 
$$
\pi_1(\mathcal{F}_{h}(\mc P^1))\cong\big<\alpha_{ij}, 1\leq i<j\leq h-1\,\big|%
(YB\,3)_{h-1},(YB\,4)_{h-1}, D_{h-1}^2=1\big>
$$
where $D_k=\alpha_{12}(\alpha_{13}\alpha_{23})(\alpha_{14}\alpha_{24}\alpha_{34}) \cdots(\alpha_{1k}\alpha_{2k}\cdots\alpha_{k-1\ k})$ 
and $(YB\,3)_n$ and $(YB\,4)_n$ are the Yang-Baxter relations (see \cite{B} and \cite{FH}): 
$$
\begin{aligned}
&(YB\,3)_n\!:\,&
\alpha_{ij}\alpha_{ik}\alpha_{jk}=\alpha_{ik}\alpha_{jk}\alpha_{ij}=
\alpha_{jk}\alpha_{ij}\alpha_{ik}, \, 1\leq i <j< k\leq n,\ \ \ \ \ \ \ \ \\
&(YB\,4)_n\!:\,&
[\alpha_{kl},\alpha_{ij}]=[\alpha_{il},\alpha_{jk}]=[\alpha_{jl},
\alpha_{jk}^{-1}\alpha_{ik}\alpha_{jk}]=[\alpha_{jl},\alpha_{kl}\alpha_{ik}
\alpha_{kl}^{-1}]=1,\\
& &  1\leq i<j<k<l\leq n .
\end{aligned}
$$
In a recent paper (\cite{BS}) Berceanu and Parveen introduced new configuration spaces.  They stratify the classical configuration spaces $\mathcal{F}_h(\mc\PP^n)$ with complex submanifolds $\mathcal{F}_{h}^{i}(\mc\PP^n)$ defined as the ordered configuration spaces of all $h$
points in $\mc\PP^n$ generating a projective subspace of dimension $i$. They prove that the fundamental groups $\pi_1(\mathcal{F}_h^i(\mc\PP^n))$ of these submanifolds are trivial except when $i=1$ providing, in this last case, a presentation similar to those of the pure braid group of the sphere.\\ 
In a subsequent paper (\cite{MPS}), authors apply similar techniques to the affine case, that is to the ordered configuration space $\mathcal{F}_h^{i,n}=\mathcal{F}_{h}^{i}(\mc^n)$ of all $h$ points in $\mc^n$ generating an affine subspace of dimension $i$. They prove that the spaces $\mathcal{F}_{h}^{i,n}$ are simply connected except for $i=1$ or $i=n=h-1$ and, in the last cases, they provide a presentation of the fundamental groups $\pi_1(\mathcal{F}_{h}^{i,n})$.\\
In this paper we generalize the result in \cite{BS} to the Grassmannian manifold $Gr(k,n)$ pa\-ra\-me\-tri\-zing $k$-dimensional subspaces of $\mc^n$. We define the $i$th ordered configuration space $\mathcal{F}_h^i(k,n)$ as the ordered configuration space of all distinct points $H_1,\ldots,H_h$ in the Grassmannian $Gr(k,n)$ such that the sum $(H_1 +\cdots+H_h)$ is an $i$-dimensional space.\\
We prove that the $i$th ordered configuration space $\mathcal{F}_h^i(k,n)$ is (when non empty) a complex sub\-ma\-ni\-fold of $Gr(k,n)^h$ and we compute its dimension.\\
As a corollary, we prove that if $n\neq hk$ and $i=$min$(n,hk)$ then the $i$th ordered configuration space $\mathcal{F}_h^i(k,n)$ has trivial fundamental group except when $n=2$, that is:
\begin{equation}
\begin{split}
&\pi_1(\mathcal{F}_h^{min(n,hk)}(k,n))=0\ \ \ {\rm if}\ (k,n)\neq(1,2)\\
&\pi_1(\mathcal{F}_1^{1}(1,2))=\pi_1(\mathcal{F}_2(\mc\PP^1)) .
\end{split}
\end{equation}
As a consequence, the fundamental group of the $i$th ordered configuration space $\mathcal{F}_h^{i}(n-1,n)$ of hyperplane arrangements of cardinality $h$ vanishes except when $n=2$.\\
Using a dual argument,  we also get that the fundamental group of the ordered configuration space of all distinct $k$-dimensional subspaces $H_1,\ldots,H_h$ in $\mc^n$ such that the intersection $(H_1 \cap \cdots \cap H_h)$ is an $i$-dimensional subspace is a simply connected manifolds when $i={\rm max}(0,n-hk)$, except when $n=2$.\\
We conjecture that similar results to that obtained in \cite{BS} for projective spaces holds also for Grassmannian manifolds and the fundamental group of the $i$th ordered configuration space $\mathcal{F}_h^i(k,n)$ vanishes except for low va\-lues of $i$. This will be the object of forthcoming publications.

\section{Main Section} 
\label{s:uno}

Let $Gr(k,n)$ be the Grassmannian manifold  pa\-ra\-me\-tri\-zing $k$-dimensional subspaces of the $n$-dimensional complex space $\mc^n$, $0<k<n$, and $\mathcal{F}_h(Gr(k,n))$ be its ordered configuration spaces. 

\paragraph{The spaces $\mathcal{F}_h^i(k,n)$.} Let's define the $i$th ordered configuration space $\mathcal{F}_h^i(k,n)$ as the space of all distinct points $H_1,\ldots,H_h$ in the Grassmannian $Gr(k,n)$ whose sum is an $i$-dimensional subspace of $\mc^n$, i.e. 
$$
\mathcal{F}_h^i(k,n)=\{(H_1,\ldots,H_h)\in\mathcal{F}_h(Gr(k,n))\mid{\rm dim}(H_1+\cdots+H_h)=i\}.
$$
It is easy to see that the following results hold:
\begin{enumerate}
\item if $h=1$ then $\mathcal{F}_1^i(k,n)$ is empty unless $i=k$, in which case $\mathcal{F}_1^k(k,n)=Gr(k,n)$;
\item if $i=1$ then $\mathcal{F}_1^i(k,n)$ is empty unless $k=h=1$ and we get $\mathcal{F}_1^1(1,n)=Gr(1,n)=\mc\PP^{n-1}$;
\item for $h\geq 2$, ${\mathcal{F}_h^i(k,n)}\neq \emptyset$ if and only if $i\geq k+1$ and $i\leq
\min(hk,n)$; 
\item for $i=hk\leq n$, then the $h$ subspaces giving a point of $\mathcal{F}_h^{hk}(k,n)$ form a direct sum;
\item for $h\geq 2$, $\mathcal{F}_h(Gr(k,n))=\mathop{\coprod}\limits_{i=2}^{n} \mathcal{F}_h^i(k,n)$;
\item for $h\geq 2$, the adjacency of the strata is given by
$$
\overline{\mathcal{F}_h^i(k,n)}=\mathcal{F}_h^i(k,n)\coprod\mathcal{F}_h^{i-1}(k,n)
\coprod\ldots\coprod\mathcal{F}_h^2(k,n).
$$
\end{enumerate}

By above remarks, it follows that the case $h=1$ is trivial, hence from now on, we will consider  $h>1$ (and hence $i>k$).\\
We want to show that $\mathcal{F}_h^i(k,n)$ is (when non empty) a complex sub\-ma\-ni\-fold of $Gr(k,n)^h$ and compute its dimension. We need to briefly recall few easy facts and introduce some notations.

\paragraph{The determinantal variety.} Let's recall that the determinantal variety $D_r(m,m')$ is the variety of $m\times m'$ matrices with complex entries of rank less than or equal to $r \leq$ min$(m,m^{\prime})$. It
is an analytic (algebraic, in fact) variety of dimension $r(m+m'-r)$ whose set of singular points is given by those matrices of rank less than $r$. From now on, $D_r(m,m')^*$ will denote the set of non-singular points of the determinantal variety $D_r(m,m')$, that is the set of $m\times m'$ matrices of rank equal to $r$.

\paragraph{A system of local coordinates for $Gr(k,n)^h$.} Let $V_0\subset\mc^n$ be a subspace of dimension dim$V_0=n-k$, then the set
$$U_{V_0}=\{H\in Gr(k,n) \mid H\oplus V_0=\mc^n\}$$
is an open dense subset of $Gr(k,n)$. \\
Let $B=\{w_1,\ldots,w_k,v_1,\ldots,v_{n-k}\}$ be a basis of $\mc^n$ such that $\{v_1,\ldots,v_{n-k}\}$ is a basis of $V_0$. We get a (complex) coordinate system on $U_{V_0}$ as follows.\\ 
Let $H$ be an element in $U_{V_0}$, then the affine subspace $V_0+w_j$ intersects $H$ in one point $u_j$ for any $j=1,\cdots,k$ and $\{u_1,\cdots,u_k\}$ form a basis of $H$. Hence $H$ is uniquely determined by a $n\times k$ matrix of the form $\begin{pmatrix}I\\A\end{pmatrix}$, where $I$ is the $k\times k$ identity matrix and $A$ is the $(n-k) \times k$ matrix of the coordinates of $u_1-w_1,\ldots,u_k-w_k$ with respect to vectors $\{v_1,\ldots,v_{n-k}\}$. The coefficients of $A$ give complex coordinates in $U_{V_0}\cong\mc^{k(n-k)}$.\\
Let $(H_1,\ldots, H_h)$ be a point in $Gr(k,n)^h$, the open sets $U_{H_1},\ldots,U_{H_h}$ in the Grassmannian manifold $Gr(n-k,n)$ have non empty intersection, that is there exists an element $V_0\in Gr(n-k,n)$ such that $V_0\oplus H_j=\mc^n$ for all $j=1,\ldots,h$. Thus, $Gr(k,n)^h$ is covered by the open sets $U_{V_0}^h$ as $V_0$ varies in $Gr(n-k,n)$. Taking a basis as defined above, each element in $U_{V_0}^h$ is uniquely determined by a $n \times hk$ matrix of the form
$\begin{pmatrix}I&I&\cdots&I\\A_1&A_2&\cdots&A_h\end{pmatrix}$ and the coefficients of $\begin{pmatrix}A_1&A_2&\cdots&A_h\end{pmatrix}$ give complex coordinates in $U_{V_0}^h\cong\mc^{hk(n-k)}$.

\paragraph{A system of local coordinates for $\mathcal{F}_h^i(k,n)$.}In terms of the above coordinates, $(H_1,\ldots, H_h)\in U_{V_0}^h$ belongs to $\mathcal{F}_h^i(k,n)$ if and only if $A_j\neq A_l$ when $j\neq l$ and rank$\begin{pmatrix}I&I&\cdots&I\\A_1&A_2&\cdots&A_h\end{pmatrix}=i$. Let us remark that 
$$\begin{aligned}
&{\rm rank}\begin{pmatrix}I&I&\cdots&I\\A_1&A_2&\cdots&A_h\end{pmatrix}&=&
{\rm rank}\begin{pmatrix}I&I&\cdots&I\\0&A_2-A_1&\cdots&A_h-A_1\end{pmatrix}\\
&&=&k+{\rm rank}\begin{pmatrix}A_2-A_1&\cdots&A_h-A_1\end{pmatrix}.
\end{aligned}$$
Then the coefficients of $B_j=A_j-A_1$ are new coordinates, in which the intersection $U_{V_0}\cap \mathcal{F}_h^i(k,n)$ corresponds, in $\mc^{hk(n-k)}$, to the product $\mc^{k(n-k)}\times D_{i-k}(n-k,hk-k)^*$  minus the closed sets given by $B_j=0$ for $2\leq j\leq h$ and by $B_j=B_l$ for $2\leq j,l\leq h$, $j\neq l$. We get the following theorem. 
\begin{theorem} The $i$th ordered configuration space $\mathcal{F}_h^i(k,n)$ is a complex submanifold of the Grassmannian manifold $Gr(k,n)$ of dimension 
\begin{equation}\label{eq:dim}
d_h^i(k,n)=i(n-i)+hk(i-k).
\end{equation}
\end{theorem}
\noindent
Equation (\ref{eq:dim}) is an easy consequence of the equality:
 $$k(n-k)+(i-k)(n-k+hk-k-(i-k))=i(n-i)+hk(i-k).$$ 
Let us remark that the dimension $d_h^i(k,n)$ attains its maximum $hk(n-k)$ if and only if $i=n$ or $i=hk$. Hence $d_h^i(k,n)$ is a strictly increasing function of $i$ when $i\leq$min$(n,hk)$. 

\paragraph{The fundamental group of $\mathcal{F}_h^{{\rm min}(n,hk)}(k,n)$.} The space  $\mathcal{F}_h^{min(n,hk)}(k,n)$ is an open subset of the ordered configuration space $\mathcal{F}_h(Gr(k,n))$ and all other (non void) $\mathcal{F}_h^j(k,n)$ have strictly lower dimension. Moreover, if $i=n$ the difference of dimensions $d_h^i(k,n)-d_h^{i-1}(k,n)$ equals $1+hk-n$ and if $i=hk$ it equals $1+n-hk$. Then if $n\neq hk$, all (non void) $\mathcal{F}_h^j(k,n)$ with $j<$min$(n,hk)$ have real codimension at least 4 in $\mathcal{F}_h(Gr(k,n))$.
Then, if $n\neq hk$ and $i=$min$(n,hk)$, the fundamental group of $\mathcal{F}_h^i(k,n)=\mathcal{F}_h(Gr(k,n))\setminus\overline{\mathcal{F}_h^{i-1}(k,n)}$
is the same as the fundamental group of $\mathcal{F}_h(Gr(k,n))$ (since, by the adjacency of the strata, the closure $\overline{\mathcal{F}_h^{i-1}(k,n)}$ is the finite union of complex subvarieties of $\mathcal{F}_h(Gr(k,n))$ of real codimension at least 4).\\
Let us recall that the complex Grassmannian manifolds $Gr(k,n)$ are simply connected and have real dimension at least 4 except $Gr(1,2)=\mc\PP^1$ and that for a simply connected manifold of real dimension at least 3 the pure braid groups vanish, i.e. $\pi_1(\mathcal{F}_h(Gr(k,n)))=0$ if $(k,n) \neq (1,2)$. We get the following corollary.
\begin{corollary} The fundamental group of the $i$th ordered configuration space $\mathcal{F}_h^{i}(k,n)$ vanishes if $n \neq hk$ and $i={\rm min}(n,hk)$ except when $n=2$ in which it is the pure braid group of the sphere.
\end{corollary}

\paragraph{The dual case.} Let $Gr(k,n)^*$ be the Grassmannian manifold parametrizing $k$-dimensional 
subspaces in the dual space $(\mc^n)^*$. Then we can define the $i$th dual ordered configuration space $\mathcal{F}_h^i(k,n)^*$ as
$$ \mathcal{F}_h^i(k,n)^*=\{(H_1,\ldots,H_h)\in \mathcal{F}_h(Gr(k,n)^*) \mid \ {\rm dim}(H_1\cap\cdots\cap H_h)=i\}.
$$
The spaces $\mathcal{F}_h^i(k,n)^*$ stratify the ordered configuration space $\mathcal{F}_h(Gr(k,n)^*)$ of the Grassmannian manifold $Gr(k,n)^*$.\\
The annihilators define homeomorphisms Ann$:Gr(n-k,n)\to Gr(k,n)^*$ which induce homeomorphisms between the $(n-i)$th ordered configuration space $F_h^{n-i}(n-k,n)$ and the $i$th dual ordered configuration space $F_h^i(k,n)^*$.\\
As a consequence the spaces $F_h^{{\rm max}(0,n-hk)}(n-k,n)^*$ are simply connected manifolds except when $n=2$. In this case the fundamental group is the pure braid group of the sphere.

\paragraph{$i$th ordered configuration spaces of hyperplane arrangements.}
If $k=n-1$ points in the ordered configuration space $\mathcal{F}_h(Gr(n-1,n))$ are $h$-uple of hyperplanes in $\mc^n$, i.e.  ordered arrangements of hyperplanes. \\
In this case, $h=1$ implies $i=n-1$ and the $i$th ordered configuration space is the Grassmannian manifold, i.e. $\mathcal{F}_1^{n-1}(n-1,n)=Gr(n-1,n)$. While $h>1$  implies $i=n$, since the sum of two (different) hyperplanes is the whole space $\mc^n$, and the following equalities hold
$$\mathcal{F}_h^{n}(n-1,n)=\mathcal{F}_h(Gr(n-1,n))=\mathcal{F}_h(\mc\PP^{n-1}).$$
Hence, the fundamental group of the $i$th ordered configuration space of hyper\-plane arrangements $\mathcal{F}_h^{i}(n-1,n)$ vanishes except when $n=2$. In this case it is the fundamental group of the sphere $\mc\PP^1$.\\
In the dual case there are homeomorphisms $\mathcal{F}_h^{i}(n-1,n)^*\cong\mathcal{F}_h^{n-i}(1,n)$ and fundamental groups of $\mathcal{F}_h^{n-i}(1,n)$ are zero except if $i=n-1$ (see \cite{BS}). Hence the space of $h$-uples of distinct hyperplanes in $\mc^n$ whose intersection has dimension equal to $i$ is simply connected except if $i=n-1$.

\bigskip

\bigskip

\textbf{Acknowledgements} The second author  gratefully acknowledge the support given to this research by the European Commission, within the 6th FP Network of Excellence ``DIME - Dynamics of Institutions and Markets in Europe'' and the Specific Targeted Research Project ``CO3 - Common Complex Collective Phenomena in Statistical Mechanics, Society, Economics and Biology''

\end{document}